\documentclass{svmult}

\usepackage{mathptmx}
\usepackage{helvet}
\usepackage{courier}
\usepackage{footmisc}
\usepackage{colonequals,amsmath,amscd,amssymb,bm}

\def\height{\operatorname{height}}

\def\Cdot{C^\bullet}
\def\Kdot{K^\bullet}
\def\bsf{{\boldsymbol{f}}}

\def\FF{\mathbb{F}}

\def\ZZ{\mathbb{Z}}
\def\calD{\mathcal{D}}
\def\ge{\geqslant}

\def\phi{\varphi}
\def\bar{\overline}

\def\del{\partial}
\def\to{\longrightarrow}
\def\mapsto{\longmapsto}

\begin{document}
\title*{A polynomial identity via differential operators}
\author{Anurag K. Singh}
\institute{A.~K.~Singh \at Department of Mathematics, University of Utah, 155 S 1400 E, Salt Lake City, UT~84112, USA \email{singh@math.utah.edu}}
\maketitle

\begin{center}
\emph{Dedicated to Professor Winfried Bruns, on the occasion of his 70th birthday}
\end{center}

\bigskip

\abstract{We give a new proof of a polynomial identity involving the minors of a matrix, that originated in the study of integer torsion in a local cohomology module.}

%%%%%%%%%%%%%%%%%%%%%%%%%%%%%%%%%%%%%%%%%%%%%%%%%%%%%%%%%%%%%%%%%%%%%%%%
\section{Introduction}
%%%%%%%%%%%%%%%%%%%%%%%%%%%%%%%%%%%%%%%%%%%%%%%%%%%%%%%%%%%%%%%%%%%%%%%%

Our study of integer torsion in local cohomology modules began in the paper~\cite{Singh:MRL}, where we constructed a local cohomology module that has $p$-torsion for each prime integer~$p$, and also studied the determinantal example $H^3_{I_2}(\ZZ[X])$ where $X$ is a $2\times 3$ matrix of indeterminates, and $I_2$ the ideal generated by its size $2$ minors. In that paper, we constructed a polynomial identity that shows that the local cohomology module $H^3_{I_2}(\ZZ[X])$ has no integer torsion; it then follows that this module is a rational vector space. Subsequently, in joint work with Lyubeznik and Walther, we showed that the same holds for all local cohomology modules of the form $H^k_{I_t}(\ZZ[X])$, where~$X$ is a matrix of indeterminates, $I_t$ the ideal generated by its size~$t$ minors, and~$k$ an integer with $k>\height I_t$, \cite[Theorem~1.2]{LSW}. In a related direction, in joint work with Bhatt, Blickle, Lyubeznik, and Zhang, we proved that the local cohomology of a polynomial ring over $\ZZ$ can have $p$-torsion for at most finitely many~$p$; we record a special case of \cite[Theorem~3.1]{BBLSZ}:

\begin{theorem}
\label{theorem:bblsz}
Let $R$ be a polynomial ring over the ring of integers, and let $f_1,\dots,f_m$ be elements of $R$. Let $n$ be a nonnegative integer. Then each prime integer that is a nonzerodivisor on the Koszul cohomology module $H^n(f_1,\dots,f_m;\,R)$ is also a nonzerodivisor on the local cohomology module $H^n_{(f_1,\dots,f_m)}(R)$.
\end{theorem}

These more general results notwithstanding, a satisfactory proof or conceptual understanding of the polynomial identity from~\cite{Singh:MRL} had previously eluded us; extensive calculations with \emph{Macaulay2} had led us to a conjectured identity, which we were then able to prove using the hypergeometric series algorithms of Petkov\v sek, Wilf, and Zeilberger~\cite{PWZ}, as implemented in \emph{Maple}. The purpose of this note is to demonstrate how techniques using differential operators underlying the papers \cite{BBLSZ} and \cite{LSW} provide the ``right'' proof of the identity, and, indeed, provide a rich source of similar identities.

We remark that there is considerable motivation for studying local cohomology of rings of polynomials with integer coefficients such as $H^k_{I_t}(\ZZ[X])$: a matrix of indeterminates~$X$ specializes to a given matrix of that size over an arbitrary commutative noetherian ring (this is where $\ZZ$ is crucial), which turns out to be useful in proving vanishing theorems for local cohomology supported at ideals of minors of arbitrary matrices. See \cite[Theorem~1.1]{LSW} for these vanishing results, that build upon the work of Bruns and Schw\"anzl~\cite{Bruns-Schwanzl}.

\section{Preliminary remarks}
\label{section:prelim}

We summarize some notation and facts. As a reference for Koszul cohomology and local cohomology, we mention \cite{Bruns-Herzog}; for more on local cohomology as a $\calD$-module, we point the reader towards \cite{Lyubeznik:Invent} and \cite{BBLSZ}.

\subsection*{Koszul and \v Cech cohomology} 

For an element $f$ in a commutative ring~$R$, the Koszul complex $\Kdot(f;\,R)$ has a natural map to the \v Cech complex $\Cdot(f;\,R)$ as follows:
\[
\CD
\Kdot(f;\,R) @. \ \colonequals\quad @. 0 @>>> R @>f>>R @>>> 0\phantom{.}\\
@VVV @. @. @| @VV\frac{1}{f}V\\ 
\Cdot(f;\,R) @. \ \colonequals\quad @. 0 @>>>R @>>> R_f @>>> 0.
\endCD
\]
For a sequence of elements $\bsf=f_1,\dots,f_m$ in $R$, one similarly obtains
\[
\CD
\Kdot(\bsf;\,R) \ \colonequals \quad \bigotimes_i K^\bullet(f_i;\,R) @>>> \bigotimes_i \Cdot(f_i;\,R) \quad \equalscolon \ C^\bullet(\bsf;\,R),
\endCD
\]
and hence, for each $n\ge0$, an induced map on cohomology modules
\begin{equation}
\label{eq:koszul:local}
\CD
H^n(\bsf;\,R) @>>> H^n_{(\bsf)}(R).
\endCD
\end{equation}

Now suppose $R$ is a polynomial ring over a field $\FF$ of characteristic $p>0$. The Frobenius endomorphism $\phi$ of $R$ induces an additive map
\[
\CD
H^n_{(\bsf)}(R) @>>> H^n_{(\bsf^p)}(R) = H^n_{(\bsf)}(R),
\endCD
\]
where $\bsf^p=f_1^p,\dots,f_m^p$. Set $R\{\phi\}$ to be the extension ring of $R$ obtained by adjoining the Frobenius operator, i.e., adjoining a generator~$\phi$ subject to the relations $\phi r=r^p\phi$ for each $r\in R$; see \cite[Section~4]{Lyubeznik:Crelle}. By an $R\{\phi\}$-module we will mean a left~$R\{\phi\}$-module. The map displayed above gives $H^n_{(\bsf)}(R)$ an $R\{\phi\}$-module structure. It is not hard to see that the image of $H^n(\bsf;\,R)$ in $H^n_{(\bsf)}(R)$ generates the latter as an $R\{\phi\}$-module; what is much more surprising is a result of \`Alvarez, Blickle, and Lyubeznik, \cite[Corollary 4.4]{ABL}, by which the image of $H^n(\bsf;\,R)$ in~$H^n_{(\bsf)}(R)$ generates the latter as a $\calD(R,\FF)$-module; see below for the definition. The result is already notable in the case~$m=1=n$, where the map~\eqref{eq:koszul:local} takes the form 
\begin{eqnarray*}
H^1(f;\,R) = R/fR & \to & R_f/R = H^1_{(f)}(R)\\
\left[1\right] & \mapsto & \left[1/f\right].
\end{eqnarray*}
By \cite{ABL}, the element $1/f$ generates $R_f$ as a $\calD(R,\FF)$-module. It is of course evident that $1/f$ generates $R_f$ as an $R\{\phi\}$-module since the elements $\phi^e(1/f)=1/f^{p^e}$ with $e\ge 0$ serve as $R$-module generators for $R_f$. See \cite{BDV} for an algorithm to explicitly construct a differential operator $\delta$ with $\delta(1/f)=1/f^{p^e}$, along with a \emph{Macaulay2} implementation.

\subsection*{Differential operators}

Let $A$ be a commutative ring, and $x$ an indeterminate; set $R=A[x]$. The divided power partial differential operator
\[
\frac{1}{k!}\frac{\del^k}{\del x^k}
\]
is the $A$-linear endomorphism of $R$ with
\[
\frac{1}{k!}\frac{\del^k}{\del x^k}(x^m) = \binom{m}{k}x^{m-k} \qquad\text{ for }\ m\ge 0,
\]
where we use the convention that the binomial coefficient $\binom{m}{k}$ vanishes if $m<k$. Note that
\[
\frac{1}{r!}\frac{\del^r}{\del x^r} \cdot \frac{1}{s!}\frac{\del^s}{\del x^s} = \binom{r+s}{r}\frac{1}{(r+s)!}\frac{\del^{r+s}}{\del x^{r+s}}.
\]

For the purposes of this paper, if $R$ is a polynomial ring over $A$ in the indeterminates $x_1,\dots,x_d$, we define the ring of $A$-linear differential operators on $R$, denoted~$\calD(R,A)$, to be the free $R$-module with basis
\[
\frac{1}{{k_1}!}\frac{\del^{k_1}}{\del x_1^{k_1}}\cdot \ \cdots\ \cdot\frac{1}{{k_d}!}\frac{\del^{k_d}}{\del x_d^{k_d}}
\qquad\text{ for }\ k_i\ge0,
\]
with the ring structure coming from composition. This is consistent with more general definitions; see \cite[16.11]{EGA4}. By a $\calD(R,A)$-module, we will mean a \emph{left} $\calD(R,A)$-module; the ring $R$ has a natural~$\calD(R,A)$-module structure, as do localizations of~$R$. For a sequence of elements $\bsf$ in $R$, the \v Cech complex $C^\bullet(\bsf;\,R)$ is a complex of~$\calD(R,A)$-modules, and hence so are its cohomology modules $H^n_{(\bsf)}(R)$. Note that for $m\ge 1$, one has
\[
\frac{1}{k!}\frac{\del^k}{\del x^k}\left(\frac{1}{x^m}\right) = (-1)^k\binom{m+k-1}{k}\frac{1}{x^{m+k}}.
\]
We also recall the Leibniz rule, which states that
\[
\frac{1}{k!}\frac{\del^k}{\del x^k}(fg) = \sum_{i+j=k}\ \frac{1}{i!}\frac{\del^i}{\del x^i}(f)\ \frac{1}{j!}\frac{\del^j}{\del x^j}(g).
\]

\section{The identity}
Let $R$ be the ring of polynomials with integer coefficients in the indeterminates
\[
\begin{pmatrix}
u & v & w\cr
x & y & z
\end{pmatrix}.
\]
The ideal $I$ generated by the size $2$ minors of the above matrix has height $2$; our interest is in proving that the local cohomology module $H^3_I(R)$ is a rational vector space. We label the minors as $\Delta_1=vz-wy$, $\Delta_2=wx-uz$, and $\Delta_3=uy-vx$. Fix a prime integer $p$, and consider the exact sequence
\[
\CD
0 @>>> R @>p>> R @>>> \bar{R} @>>> 0,
\endCD
\]
where $\bar{R}=R/pR$. This induces an exact sequence of local cohomology modules
\[
\CD
@>>> H^2_I(R) @>\pi>> H^2_I(\bar{R}) @>>> H^3_I(R) @>p>> H^3_I(R) @>>> H^3_I(\bar{R}) @>>> 0.
\endCD
\]
The ring $\bar{R}/I\bar{R}$ is Cohen-Macaulay of dimension $4$, so \cite[Proposition~III.4.1]{PS} implies that~$H^3_I(\bar{R})=0$. As $p$ is arbitrary, it follows that $H^3_I(R)$ is a divisible abelian group. To prove that it is a rational vector space, one needs to show that multiplication by $p$ on~$H^3_I(R)$ is injective, equivalently that $\pi$ is surjective. We first prove this using the identity~\eqref{eq:identity} below, and then proceed with the proof of the identity.

For each $k\ge 0$, one has
\begin{multline}
\label{eq:identity}
\sum_{i,j\ge0}\binom{k}{i+j}\binom{k+i}{k}\binom{k+j}{k}\frac{(-wx)^i(vx)^ju^{k+1}}{\Delta_2^{k+1+i}\Delta_3^{k+1+j}}\\
+
\sum_{i,j\ge0}\binom{k}{i+j}\binom{k+i}{k}\binom{k+j}{k}\frac{(-uy)^i(wy)^jv^{k+1}}{\Delta_3^{k+1+i}\Delta_1^{k+1+j}}\\
+
\sum_{i,j\ge0}\binom{k}{i+j}\binom{k+i}{k}\binom{k+j}{k}\frac{(-vz)^i(uz)^jw^{k+1}}{\Delta_1^{k+1+i}\Delta_2^{k+1+j}}
\ \ =\ \ 0.
\end{multline}
Since the binomial coefficient $\binom{k}{i+j}$ vanishes if $i$ or $j$ exceeds $k$, this equation may be rewritten as an identity in the polynomial ring $\ZZ[u,v,w,x,y,z]$ after multiplying by $(\Delta_1\Delta_2\Delta_3)^{2k+1}$.

Computing $H^2_I(R)$ as the cohomology of the \v Cech complex $\Cdot(\Delta_1,\Delta_2,\Delta_3;\,R)$, equation~\eqref{eq:identity} gives a $2$-cocycle in
\[
C^2(\Delta_1,\Delta_2,\Delta_3;\,R)=R_{\Delta_1\Delta_2}\oplus R_{\Delta_1\Delta_3} \oplus R_{\Delta_2\Delta_3};
\]
we denote the cohomology class of this cocycle in $H^2_I(R)$ by $\eta_k$. When $k=p^e-1$, one has
\[
\binom{k}{i+j}\binom{k+i}{k}\binom{k+j}{k}\equiv0\mod p\quad\text{ for }\ (i,j)\neq(0,0),
\]
so~\eqref{eq:identity} reduces modulo $p$ to
\[
\frac{u^{p^e}}{\Delta_2^{p^e}\Delta_3^{p^e}} + \frac{v^{p^e}}{\Delta_3^{p^e}\Delta_1^{p^e}} + \frac{w^{p^e}}{\Delta_1^{p^e}\Delta_2^{p^e}} \equiv0\mod p,
\]
and the cohomology class $\eta_{p^e-1}$ has image
\[
\pi(\eta_{p^e-1})=\left[\left(\frac{w^{p^e}}{\Delta_1^{p^e}\Delta_2^{p^e}}, \frac{-v^{p^e}}{\Delta_1^{p^e}\Delta_3^{p^e}}, \frac{u^{p^e}}{\Delta_2^{p^e}\Delta_3^{p^e}}\right)\right] \qquad\text{ in }\ H^2_I(\bar{R}).
\]
Since $\bar{R}$ is a regular ring of positive characteristic, $H^2_I(\bar{R})$ is generated as an $\bar{R}\{\phi\}$-module by the image of
\[
\CD
H^2(\Delta_1,\Delta_2,\Delta_3;\,\bar{R})@>>>H^2_I(\bar{R}).
\endCD
\]
The Koszul cohomology module $H^2(\Delta_1,\Delta_2,\Delta_3;\,\bar{R})$ is readily seen to be generated, as an~$\bar{R}$-module, by elements corresponding to the relations
\[
u\Delta_1+v\Delta_2+w\Delta_3=0\qquad\text{ and }\qquad x\Delta_1+y\Delta_2+z\Delta_3=0.
\]
These two generators of $H^2(\Delta_1,\Delta_2,\Delta_3;\,\bar{R})$ map, respectively, to
\[
\alpha:=\left[\left(\frac{w}{\Delta_1\Delta_2}, \frac{-v}{\Delta_1\Delta_3}, \frac{u}{\Delta_2\Delta_3}\right)\right]
\quad\text{ and }\quad
\beta:=\left[\left(\frac{z}{\Delta_1\Delta_2}, \frac{-y}{\Delta_1\Delta_3}, \frac{x}{\Delta_2\Delta_3}\right)\right]
\]
in $H^2_I(\bar{R})$. Thus, $H^2_I(\bar{R})$ is generated over $\bar{R}$ by $\phi^e(\alpha)$ and $\phi^e(\beta)$ for $e\ge0$. But
\[
\phi^e(\alpha) = \pi(\eta_{p^e-1})
\]
is in the image of $\pi$, and hence so is $\phi^e(\beta)$ by symmetry. Thus, $\pi$ is surjective.

\subsection*{The proof of the identity}

We start by observing that $C^2(\Delta_1,\Delta_2,\Delta_3;\,R)$ is a $\calD(R,\ZZ)$-module. The element
\[
\left(\frac{w}{\Delta_1\Delta_2}, \frac{-v}{\Delta_1\Delta_3}, \frac{u}{\Delta_2\Delta_3}\right)
\]
is a $2$-cocycle in $C^2(\Delta_1,\Delta_2,\Delta_3;\,R)$ since 
\begin{equation}
\label{eq:det}
\frac{w}{\Delta_1\Delta_2}+\frac{v}{\Delta_1\Delta_3}+\frac{u}{\Delta_2\Delta_3}=0.
\end{equation}
We claim that the identity~\eqref{eq:identity} is simply the differential operator 
\[
D=\frac{1}{k!}\frac{\del^k}{\del u^k} \cdot \frac{1}{k!}\frac{\del^k}{\del y^k} \cdot \frac{1}{k!}\frac{\del^k}{\del z^k}
\]
applied termwise to~\eqref{eq:det}; we first explain the choice of this operator: set $k=p^e-1$, and consider $\bar{D}=D\mod p$ as an element of
\[
\calD(R,\ZZ)/p\calD(R,\ZZ) = \calD(R/pR,\ZZ/p\ZZ).
\]
It is an elementary verification that
\begin{eqnarray*}
\bar{D}(u\Delta_2^{p^e-1}\Delta_3^{p^e-1}) \ \equiv\ u^{p^e}&\\
\bar{D}(v\Delta_3^{p^e-1}\Delta_1^{p^e-1}) \ \equiv\ v^{p^e}& \mod p.\\
\bar{D}(w\Delta_1^{p^e-1}\Delta_2^{p^e-1}) \ \equiv\ w^{p^e}&
\end{eqnarray*}
Since $k<p^e$, the differential operator $\bar{D}$ is $\bar{R}^{p^e}$-linear; dividing the above equations by $\Delta_2^{p^e}\Delta_3^{p^e}$, $\Delta_3^{p^e}\Delta_1^{p^e}$, and $\Delta_1^{p^e}\Delta_2^{p^e}$ respectively, we obtain
\[
\bar{D}\left(\frac{w}{\Delta_1\Delta_2}, \frac{-v}{\Delta_1\Delta_3}, \frac{u}{\Delta_2\Delta_3}\right)
\equiv
\left(\frac{w^{p^e}}{\Delta_1^{p^e}\Delta_2^{p^e}}, \frac{-v^{p^e}}{\Delta_1^{p^e}\Delta_3^{p^e}}, \frac{u^{p^e}}{\Delta_2^{p^e}\Delta_3^{p^e}}\right) \mod p,
\]
which maps to the desired cohomology class $\phi^e(\alpha)$ in $H^2_I(\bar{R})$. Of course, the operator $D$ is not unique in this regard.

Using elementary properties of differential operators recorded in \S\ref{section:prelim}, we have
\begin{align*}
D\left(\frac{v}{\Delta_3\Delta_1}\right)
&= \frac{1}{k!}\frac{\del^k}{\del u^k} \cdot \frac{1}{k!}\frac{\del^k}{\del y^k} \cdot \frac{1}{k!}\frac{\del^k}{\del z^k} 
\ \left[\frac{v}{(uy-vx)(vz-wy)}\right]\\
&= \frac{1}{k!}\frac{\del^k}{\del u^k} \cdot \frac{1}{k!}\frac{\del^k}{\del y^k}
\ \left[\frac{v(-v)^k}{(uy-vx)(vz-wy)^{k+1}}\right]\\
&= \frac{1}{k!}\frac{\del^k}{\del y^k}
\ \left[\frac{v(-v)^k(-y)^k}{(uy-vx)^{k+1}(vz-wy)^{k+1}}\right]\\
&= v^{k+1} \frac{1}{k!}\frac{\del^k}{\del y^k} \ \left[\frac{y^k}{(uy-vx)^{k+1}(vz-wy)^{k+1}}\right]\\
&= v^{k+1}\sum_{i,j} 
\left[\frac{1}{i!}\frac{\del^i}{\del y^i}\frac{1}{(uy-vx)^{k+1}}\right]
\left[\frac{1}{j!}\frac{\del^j}{\del y^j}\frac{1}{(vz-wy)^{k+1}}\right]
\left[\frac{1}{(k-i-j)!}\frac{\del^{k-i-j}}{\del y^{k-i-j}}y^k\right]\\
&= v^{k+1}\sum_{i,j} 
\binom{k+i}{i}\frac{(-u)^i}{(uy-vx)^{k+1+i}}
\binom{k+j}{j}\frac{w^j}{(vz-wy)^{k+1+j}}
\binom{k}{i+j}y^{i+j}\\
&= v^{k+1}\sum_{i,j} \binom{k+i}{i}\binom{k+j}{j}\binom{k}{i+j}
\frac{(-uy)^i(wy)^j}{\Delta_3^{k+1+i}\Delta_1^{k+1+j}}.
\end{align*}

A similar calculation shows that
\[
D\left(\frac{w}{\Delta_1\Delta_2}\right)=w^{k+1}\sum_{i,j} \binom{k+i}{i}\binom{k+j}{j}\binom{k}{i+j}
\frac{(-vz)^i(uz)^j}{\Delta_1^{k+1+i}\Delta_2^{k+1+j}}.
\]
It remains to evaluate $\displaystyle{D\left(\frac{u}{\Delta_2\Delta_3}\right)}$; we reduce this to the previous calculation as follows. First note that the differential operators $\displaystyle{\frac{\del}{\del u}\cdot\frac{\del}{\del y}}$ and $\displaystyle{\frac{\del}{\del v}\cdot\frac{\del}{\del x}}$ commute; it is readily checked that they agree on $\displaystyle{\frac{u}{\Delta_2\Delta_3}}$. Consequently the operators
\[
\frac{1}{k!}\frac{\del^k}{\del u^k} \cdot \frac{1}{k!}\frac{\del^k}{\del y^k} \cdot \frac{1}{k!}\frac{\del^k}{\del z^k}
\quad\text{ and }\quad
\frac{1}{k!}\frac{\del^k}{\del v^k} \cdot \frac{1}{k!}\frac{\del^k}{\del z^k} \cdot \frac{1}{k!}\frac{\del^k}{\del x^k}
\]
agree on $\displaystyle{\frac{u}{\Delta_2\Delta_3}}$ as well. But then
\[
D\left(\frac{u}{\Delta_2\Delta_3}\right)
=\frac{1}{k!}\frac{\del^k}{\del v^k} \cdot \frac{1}{k!}\frac{\del^k}{\del z^k} \cdot \frac{1}{k!}\frac{\del^k}{\del x^k}
\ \left[\frac{u}{(wx-uz)(uy-vx)}\right]
\]
which, using the previous calculation and symmetry, equals 
\[
u^{k+1}\sum_{i,j} \binom{k+i}{i}\binom{k+j}{j}\binom{k}{i+j}
\frac{(-wx)^i(vx)^j}{\Delta_2^{k+1+i}\Delta_3^{k+1+j}}.
\]

\subsection*{Identities in general}

Suppose $\bsf=f_1,\dots,f_m$ are elements of a polynomial ring $R$ over~$\ZZ$, and $g_1,\dots,g_m$ are elements of $R$ such that
\[
g_1f_1+\cdots+g_mf_m=0.
\]
Then, for each prime integer $p$ and $e\ge0$, the Frobenius map on $\bar{R}=R/pR$ gives
\begin{equation}
\label{eq:frob}
g_1^{p^e}f_1^{p^e}+\cdots+g_m^{p^e}f_m^{p^e}\equiv 0\mod p.
\end{equation}
Now suppose $p$ is a nonzerodivisor on the Koszul cohomology module~$H^m(\bsf;\,R)$. Then Theorem~\ref{theorem:bblsz} implies that~\eqref{eq:frob} \emph{lifts} to an equation
\begin{equation}
\label{eq:lift}
G_1f_1^N+\cdots+G_mf_m^N=0
\end{equation}
in $R$ in the sense that the cohomology class corresponding to~\eqref{eq:lift} in $H^{m-1}_{(\bsf)}(R)$ maps to the cohomology class corresponding to~\eqref{eq:frob} in $H^{m-1}_{(\bsf)}(\bar{R})$.

\begin{acknowledgement}
NSF support under grant DMS~1500613 is gratefully acknowledged. This paper owes an obvious intellectual debt to our collaborations with Bhargav Bhatt, Manuel Blickle, Gennady Lyubeznik, Uli Walther, and Wenliang Zhang; we take this opportunity to thank our coauthors.
\end{acknowledgement}

%%%%%%%%%%%%%%%%%%%%%%%%%%%%%%%%%%%%%%%%%%%%%%%%%%%%%%%%%%%%%%%%%%%%%%%%

\end{document}